\documentclass[12pt]{amsart}
\usepackage{amscd,amssymb,amsthm,verbatim}
\def\bp{\star} 

\def\barint{\mathop{-\mkern-19.5mu\int}}
\def\av#1{{\langle #1 \rangle}} 

\newcommand{\const}{\operatorname{const.}}

\newcommand{\diam}{\operatorname{diam}}

\newcommand{\dvol}{\operatorname{dvol}}

\newcommand{\Hess}{\operatorname{Hess}}

\newcommand{\Id}{\operatorname{Id}}

\newcommand{\Lip}{\operatorname{Lip}}

\newcommand{\R}{{\mathbb R}}

\newcommand{\Ric}{\operatorname{Ric}}

\newcommand{\supp}{\operatorname{supp}}

\newcommand{\Tr}{\operatorname{Tr}}

\newcommand{\vol}{\operatorname{vol}}
\newcommand{\Z}{{\mathbb Z}}

\newcommand{\DC}{{\cal D}{\cal C}}
\numberwithin{equation}{section}
\theoremstyle{plain}
\newtheorem{definition}[equation]{Definition}

\newtheorem{theorem}[equation]{Theorem}
\newtheorem{proposition}[equation]{Proposition}
\newtheorem{corollary}[equation]{Corollary}

\theoremstyle{remark}

\newtheorem{remark}[equation]{Remark}
\newtheorem{example}[equation]{Example}
\newtheorem{wishlist}[equation]{Wishlist}

\usepackage{graphicx}
\input{amssym.def}
\input{amssym}
\input{epsf}
\usepackage{a4wide}

\def\R{\mathbb R}

\def\Z{\mathbb Z}

\def\cal{\mathcal}

\def\Hess{\mathop{{\rm Hess}\,}}

\def\Id{{\rm Id}\,}

\def\2dr#1#2{\left. \frac{d^2}{d{#1}^2} \right |_{#2}}
\def\d2#1{\frac{d^2}{d{#1}^2}}

\DeclareMathOperator*{\grad}{grad}

\def\begeq{\begin{equation}}
\def\endeq{\end{equation}}
\def\begar{\begin{eqnarray}}
\def\endar{\end{eqnarray}}
\def\begar*{\begin{eqnarray*}}
\def\endar*{\end{eqnarray*}}
\def\begal{\begin{align}}
\def\endal{\end{align}}
\def\begal*{\begin{align*}}
\def\endal*{\end{align*}}

\theoremstyle{definition}

\theoremstyle{remark}

\newtheorem*{Thm*}{Theorem}
\newtheorem*{Lem*}{Lemma}
\newtheorem*{Conj*}{Conjecture}
\newtheorem*{Cor*}{Corollary}
\newtheorem*{Def*}{Definition}
\newtheorem*{Prop*}{Proposition}
\newtheorem*{Exo*}{Exercise}
\newtheorem*{Exs*}{Examples}
\newtheorem*{Ex*}{Example}
\newtheorem*{Rk*}{Remark}
\newtheorem*{Rks*}{Remarks}

\begin{document}

\title{Optimal transport and Ricci curvature for
metric-measure spaces} 

\author{John Lott}
\address{Department of Mathematics\\
University of Michigan\\
Ann Arbor, MI  48109-1109\\
USA} \email{lott@umich.edu}

\thanks{The author was partially
supported by NSF grant DMS-0604829 during the
writing of this article}
\date{December 15, 2006}

\begin{abstract} We survey work of Lott-Villani and
Sturm on lower Ricci curvature
bounds for metric-measure spaces.
\end{abstract}

\maketitle
An intriguing question is whether one
can extend notions of smooth Riemannian geometry to
general metric spaces.  Besides the inherent interest,
such extensions sometimes allow one to prove
results about smooth Riemannian manifolds, using
compactness theorems.

There is a good notion of a metric space having
``sectional curvature bounded below by $K$'' or
``sectional curvature bounded above by $K$'', due to
Alexandrov.  We refer to the articles of
Petrunin and Buyalo-Schroeder in this
volume for further information on these two topics.
In this article we address the issue of whether there
is a good notion of a metric space having
``Ricci curvature bounded below by $K$''.

A motivation for this question comes from Gromov's
precompactness theorem \cite[Theorem 5.3]{Gromov (1999)}. 
Let ${\mathcal M}$ denote
the set of compact metric spaces (modulo isometry)
with the Gromov-Hausdorff topology. The precompactness
theorem says that given
$N \in \Z^+$, $D < \infty$ and $K \in \R$, the
subset of ${\mathcal M}$ consisting of
closed Riemannian manifolds $(M, g)$ with
$\dim(M) = N$, $\Ric \ge Kg$ and
$\diam \le D$, is precompact.
The limit points in ${\mathcal M}$ of this subset
will be metric spaces of Hausdorff dimension
at most $N$, but generally are not manifolds. However,
one would like to say that 
in some generalized sense
they do have Ricci curvature bounded 
below by $K$. Deep
results about the structure of such limit
points, which we call Ricci limits, were
obtained by Cheeger and Colding
\cite{Cheeger-Colding (1997),Cheeger-Colding II (2000),Cheeger-Colding (2000)}.
We refer to the article
of Guofang Wei in this volume for further information.

In the work of Cheeger and Colding, and in 
earlier work of Fukaya \cite{Fukaya (1987)}, it turned out to be useful to
consider not just metric spaces, but rather metric
spaces equipped with measures. Given a compact
metric space $(X, d)$, let $P(X)$ denote the
set of Borel probability measures on $X$. That is,
$\nu \in P(X)$ means that $\nu$ is a nonnegative
Borel measure on $X$ with $\int_X d\nu \: = \: 1$.
We put the weak-$*$ topology on $P(X)$, so
$\lim_{i \rightarrow \infty} \nu_i \: = \: \nu$ if and only if
for all $f \in C(X)$, we have
$\lim_{i \rightarrow \infty} \int_X f \: d\nu_i \: = \: 
\int_X f \: d\nu$. Then $P(X)$ is compact.

\begin{definition} \label{def1}
A compact {\em metric measure space} is a triple
$(X, d, \nu)$ where $(X,d)$ is a compact metric space and $\nu \in P(X)$.
\end{definition}
 
\begin{definition} \label{def2}
Given two compact metric spaces $(X_1,d_1)$ and $(X_2,d_2)$, 
an $\epsilon$-Gromov-Hausdorff 
approximation from $X_1$ to $X_2$ is a (not necessarily continuous)
map $f \: : \: X_1 \rightarrow X_2$ so that 

(i)  For all $x_1, x_1^\prime \in X_1$, 
$\bigl|d_{2}(f(x_1), f(x_1^\prime)) \: - \: d_{1}(x_1, x_1^\prime)\bigr|
\: \le \: \epsilon$. 

(ii) For all $x_2 \in X_2$, there is an $x_1 \in X_1$ so that
$d_{2}(f(x_1), x_2) \: \le \: \epsilon$.

A sequence
$\{(X_i, d_i, \nu_i)\}_{i=1}^\infty$ of compact
metric-measure spaces converges to $(X,d, \nu)$ in the
{\em measured Gromov-Hausdorff topology} 
if there are Borel $\epsilon_i$-approximations $f_i:X_i\to X$, with
$\lim_{i \rightarrow \infty} \epsilon_i \: = \: 0$, so that
$\lim_{i \rightarrow \infty} (f_i)_* \nu_i \: = \: \nu$ in $P(X)$.
\end{definition}

\begin{remark} \label{rmk1}
There are other interesting topologies on the set of
metric-measure spaces, discussed in \cite[Chapter $3 \frac12$]{Gromov (1999)}.
\end{remark}

If $M$ is a compact manifold with Riemannian metric $g$ then we also
let $(M, g)$ denote the underlying metric space. There is
a canonical probability
measure on $M$ given by the normalized volume form
$\frac{\dvol_M}{\vol(M)}$. 
One can easily extend Gromov's
precompactness theorem to say that
given
$N \in \Z^+$, $D < \infty$ and $K \in \R$, the
triples $\left( M, g, \frac{\dvol_M}{\vol(M)} \right)$ with
$\dim(M) = N$, $\Ric \ge Kg$ and
$\diam \le D$ form a precompact subset in the measured
Gromov-Hausdorff (MGH) topology.  The limit points of this subset
are now metric-measure spaces $(X, d, \nu)$. One would
like to say that they have 
``Ricci curvature bounded below by $K$'' in some
generalized sense.

The metric space $(X,d)$ of a Ricci limit is necessarily
a length space. Hereafter we mostly restrict our attention to
length spaces.
So the question that we address is whether there is a good notion
of a compact measured length space $(X, d, \nu)$ having
``Ricci curvature bounded below by $K$''.  The word
``good'' is a bit ambiguous here, but we would 
like our definition to have the following properties. 
\begin{wishlist} \label{wishlist}
1. If $\{(X_i, d_i, \nu_i)\}_{i=1}^\infty$ is a sequence of compact
measured length spaces with ``Ricci curvature bounded below by $K$''
and $\lim_{i \rightarrow \infty} (X_i, d_i, \nu_i) \: = \: 
(X, d, \nu)$ in the measured Gromov-Hausdorff topology then
$(X,d,\nu)$ has ``Ricci curvature bounded below by $K$''. \\
2. If $(M, g)$ is a compact Riemannian manifold
then the triple
$\left( M, g, \frac{\dvol_M}{\vol(M)} \right)$ has 
``Ricci curvature bounded below by $K$'' if and only if
$\Ric \ge Kg$ in the usual sense. \\
3. One can prove some nontrivial results about
measured length spaces having
``Ricci curvature bounded below by $K$''. 
\end{wishlist}

It is not so easy to come up with a definition that satisfies all
of these properties.  One possibility would be to say that
$(X, d, \nu)$ has ``Ricci curvature bounded below by $K$''
if and only if it is an MGH limit of Riemannian manifolds with
$\Ric \ge K g$, but this is a bit
tautological. We want instead a definition that
depends in an intrinsic way on $(X,d,\nu)$. We refer to
\cite[Appendix 2]{Cheeger-Colding (1997)} 
for further discussion of the problem.

In fact, it will turn out that we will want to specify an
effective dimension $N$, possibly infinite, of the 
measured length space.
That is, we want to define a notion of $(X, d, \nu)$ having
``$N$-Ricci curvature bounded below by $K$'', where
$N$ is a parameter that is part of the definition.
The need to input the parameter $N$ can be seen from the Bishop-Gromov
inequality for complete
$n$-dimensional Riemannian manifolds with nonnegative Ricci
curvature. It says that $r^{-n} \vol(B_r(m))$ is nonincreasing in $r$,
where $B_r(m)$ is the $r$-ball centered at $m$.
We will want a Bishop-Gromov-type inequality to hold in the
length space setting, but
when we go from manifolds to 
length spaces there is no {\it a priori} value for the parameter $n$.
Hence for each $N \in [1, \infty]$, there will be a notion of 
$(X, d, \nu)$ having
``$N$-Ricci curvature bounded below by $K$''.

The goal now is to find some property which we know holds for
$N$-dimensional
Riemannian manifolds with Ricci curvature bounded below, and turn it
into a definition for measured length spaces.  A geometer's
first inclination may be to just use the Bishop-Gromov
inequality, at least if $N < \infty$, for example to say that $(X, d, \nu)$ has
``nonnegative $N$-Ricci curvature'' if and only if
for each $x \in \supp(\nu)$, $r^{-N} \: \nu(B_r(x))$ is
nonincreasing in $r$. Although this is the simplest possibility,
it turns out that it is not satisfactory; see Remark \ref{rmk8}.
Instead, we will derive a Bishop-Gromov inequality as
part of a more subtle definition.

The definition that we give in this paper may seem to come from left
field, at least from the viewpoint of standard geometry.  It comes from
a branch of applied mathematics called {\em optimal transport},
which can be informally considered to be the study of moving dirt
around.  The problem originated with Monge in the paper
\cite{Monge}, whose title translates into English as
``On the theory of excavations and fillings''.  (In that paper
Monge also introduced the idea of a line of curvature of a surface.)
The problem that Monge
considered was how to transport a ``before'' dirtpile to
an ``after'' dirtpile with minimal total ``cost'', 
where he took the cost of
transporting a unit mass of dirt between points $x$ and $y$ to be $d(x,y)$. 
Such a transport $F : X \rightarrow X$
is called a {\em Monge transport}.
An account of Monge's life, and his unfortunate political choices,
is in \cite{Bell}.

Since Monge's time, there has been considerable work on
optimal transport.  Of course, the original case of interest was optimal
transport on Euclidean space.  Kantorovich introduced a important
relaxation of Monge's original problem, in which not all of the
dirt from a given point $x$ has to go to a single point $y$.  That is,
the dirt from $x$ is allowed to be spread out over the space.
Kantorovich showed that there is always an optimal transport scheme in
his sense. (Kantorovich won a 1975 Nobel Prize in economics.)
We refer to the book \cite{TOT} for a lively and detailed
account of optimal transport. 

In Section \ref{sec1} we summarize some optimal transport results from a
modern perspective.  We take the cost function of transporting
a unit mass of dirt to be $d(x,y)^2$ instead of Monge's
cost function $d(x,y)$. The relation to Ricci curvature comes from
work of Otto-Villani \cite{Otto-Villani (2000)} and 
Cordero-Erausquin-McCann-Schmuckenschl\"ager \cite{CEMS01}.
They showed that optimal transport on a
Riemannian manifold is affected by the Ricci tensor.
To be a bit more precise, the Ricci curvature affects the
convexity of certain entropy functionals along an optimal
transport path.  Details are in Section \ref{sec2}.

The idea now, implemented independently by Lott-Villani and
Sturm, is to {\em define} the property ``$N$-Ricci curvature bounded below by $K$'',
for a measured length space $(X,d,\nu)$,
in terms of the convexity of certain entropy functionals
along optimal transport paths in the auxiliary space $P(X)$. We present the
definition and its initial properties in Section \ref{sec3}.
We restrict in that section to
the case $K=0$, where the discussion becomes a bit simpler.
We show that Condition 1. from 
Wishlist \ref{wishlist} is satisfied.  In Section \ref{sec4} we show
that Condition 2. from Wishlist \ref{wishlist} is satisfied.
In Section \ref{sec5} we give the definition of $(X,d,\nu)$ having
$N$-Ricci curvature bounded below by $K$, for $K \in \R$.

Concerning Condition 3. of the
Wishlist, in Sections \ref{sec3}, \ref{sec4} and \ref{sec5}
we give some geometric results that one can prove about measured
length spaces with Ricci curvature bounded below. In
particular, there are applications to Ricci limit spaces.
In Section \ref{sec6} we give some analytic results.
In Section \ref{sec7} we discuss some further issues.

We mostly focus on results
from \cite{LV} and \cite{LV2}, mainly because of the author's
familiarity with those papers.  However, we emphasize that
many parallel results were obtained independently by 
Karl-Theodor Sturm in \cite{Sturm2,Sturm3}. Background information on optimal
transport is in \cite{TOT} and \cite{TOT2}. The latter book also
contains a more detailed exposition of
some of the topics of this survey.

I thank C\'edric Villani for an enjoyable collaboration.

\section{Optimal transport} \label{sec1}

Let us state the Kantorovich transport problem.
We take $(X, d)$ to be a compact metric space. Our
``before'' and ``after'' dirtpiles are measures
$\mu_0, \mu_1 \in P(X)$. They both have mass one.
We want to move the total amount of dirt from $\mu_0$ to $\mu_1$
most efficiently. A moving scheme, maybe not optimal, will be 
called a {\em transference plan}. Intuitively, it amounts to specifying how
much dirt is moved from a point $x_0$ to a point $x_1$. That is, we have a
probability measure
$\pi \in P(X\times X)$, which we informally write as $\pi(x_0,x_1)$.
The statement
that $\pi$ does indeed transport $\mu_0$ to $\mu_1$ translates to
the condition that 
\begin{equation} \label{eqn1}
(p_0)_* \pi \: = \: \mu_0,\qquad (p_1)_* \pi \: = \: \mu_1,
\end{equation}
where $p_0, p_1 \: : \: X \times X \rightarrow X$ are projections
onto the first and second factors, respectively.

We will use optimal transport with quadratic cost function
(square of the distance). The total cost of the transference plan
$\pi$ is given by adding the contributions of $d(x_0, x_1)^2$ with
respect to $\pi$. 
Taking the infimum of this with respect to $\pi$ 
gives the square of the {\em Wasserstein distance} 
$W_2(\mu_0, \mu_1)$
between $\mu_0$ and $\mu_1$, i.e.
\begin{equation} \label{var}
W_2(\mu_0, \mu_1)^2 \: = \:
\inf_{\pi}  \: \int_{X \times X} d(x_0, x_1)^2 \: d\pi(x_0, x_1),
\end{equation}
where $\pi$ ranges over the set of all transference plans between
$\mu_0$ and $\mu_1$. Any minimizer $\pi$ for this variational problem
is called an {\em optimal transference plan}.

In \eqref{var}, one can replace the infimum by the minimum
\cite[Proposition 2.1]{TOT}, i.e. there always exists
(at least) one optimal transference plan.
It turns out that $W_2$ is a metric on $P(X)$. 
The topology that it induces on $P(X)$ is the
weak-$*$ topology~\cite[Theorems 7.3 and 7.12]{TOT}.
When equipped with the metric $W_2$, $P(X)$ is a compact metric space.
In this way, to each compact metric space $X$ we have assigned
another compact metric space $P(X)$.
The {\em Wasserstein space} $(P(X), W_2)$ seems to be a very natural 
object in mathematics. It generally has infinite topological or
Hausdorff dimension.  (If $X$ is a finite set then $P(X)$ is a
simplex, with a certain metric.) It is always contractible, as can
be seen by fixing a measure $\mu_0 \in P(X)$ and 
linearly contracting
other measures $\mu \in P(X)$ to $\mu_0$ by
$t \rightarrow t\mu_0 \: + \: (1-t) \mu$.

\begin{proposition} \cite[Corollary 4.3]{LV} \label{prop1}
If $\lim_{i \rightarrow \infty} (X_i, d_i) \: = \: 
(X,d)$ in the Gromov-Hausdorff topology 
then $\lim_{i \rightarrow \infty} (P(X_i), W_2) \: = \: 
(P(X), W_2)$ in the Gromov-Hausdorff topology.
\end{proposition}

A {\em Monge transport} is a transference plan coming from a map
$F \: : \: X \rightarrow X$ with $F_* \mu_0 \: = \: \mu_1$, given
by $\pi \: = \: (\Id, F)_* \mu_0$. In general an optimal transference
plan does not have to be a Monge transport, although this may
be true under some assumptions.

What does optimal transport look like in Euclidean space $\R^n$?
Suppose that $\mu_0$ and $\mu_1$ are compactly supported and
absolutely continuous with respect to Lebesgue measure.  Brenier 
\cite{Brenier} and
Rachev-R\"uschendorf \cite{RR} showed that
there is a unique optimal transference plan between $\mu_0$ and
$\mu_1$, which is a Monge
transport.  Furthermore, there is a convex function $V$ on 
$\R^n$ so that for almost all $\vec{x}$, the Monge transport is given by
$F(\vec{x}) \: = \: \vec{\nabla}_{\vec{x}} V$. So to find the optimal
transport, one finds a convex function $V$ such that
the pushforward, under the map
$\nabla V \: : \: \R^n \rightarrow \R^n$, sends $\mu_0$ to
$\mu_1$. This solves the Monge problem for such measures,
under our assumption of quadratic cost function.  The solution to
the original problem of Monge, with linear cost function, is
more difficult; see \cite{Evans-Gangbo}.

The statement of the Brenier-Rachev-R\"uschendorf theorem may sound like
anathema to a geometer.  One is identifying the gradient
of $V$ (at $\vec{x}$), which is a vector, with the image of 
$\vec{x}$ under a map, which is a
point. Because of this, it is not evident how to extend even the
statement of the
theorem if one wants to do optimal transport on 
a Riemannian manifold. The extension was done by McCann \cite{McCann (2001)}.
The key point is that on $\R^n$, we can write
$\vec{\nabla}_{\vec{x}} V \: = \: \vec{x} - \vec{\nabla}_{\vec{x}} \phi$,
where $\phi(\vec{x}) \: = \: \frac{|\vec{x}|^2}{2} \: - \: V(\vec{x})$.
To understand the relation between $V$ and $\phi$, we note that if the
convex function $V$ were smooth then $\phi$ would have Hessian
bounded above by the identity.  On a Riemannian manifold $(M, g)$,
McCann's theorem says that
an optimal transference plan between two compactly supported
absolutely continuous
measures is a Monge
transport $F$ that satisfies $F(m) \: = \: \exp_m (- \nabla_m \phi)$
for almost all $m$,
where $\phi$ is a function on $M$ with Hessian bounded above by $g$
in a generalized sense. More precisely, $\phi$
is $\frac{d^2}{2}$-concave in the sense that it can
be written in the form
\begin{equation} \label{eqn2}
\phi(m) \: = \: \inf_{m^\prime \in M} 
\left( \frac{d(m, m^\prime)^2}{2} \: - \: \widetilde{\phi}(m^\prime) \right) 
\end{equation}
for some function $\widetilde{\phi} \: : \: M \rightarrow [-\infty, \infty)$.

Returning to the metric space setting, if $(X, d)$ is a compact
length space and one has an optimal transference plan
$\pi$ then one would physically perform the transport
by picking up pieces of dirt in $X$ and moving them along minimal
geodesics to other points in $X$, in a way consistent with 
the transference plan
$\pi$. The transference plan $\pi$ tells us how much dirt has to go
from $x_0$ to $x_1$, but does not say anything about which
minimal geodesics from $x_0$ to $x_1$ we should actually use. After making
such a choice of minimizing geodesics, we obtain a $1$-parameter family of
measures $\{\mu_t\}_{t \in [0,1]}$ by stopping the
physical transport procedure at time $t$ and looking
at where the dirt is.
This suggests looking at $(P(X), W_2)$ as a length space.

\begin{proposition} 
\cite[Corollary 2.7]{LV},\cite[Proposition 2.10(iii)]{Sturm2} \label{prop1.5}
If $(X, d)$ is a compact length space then $(P(X), W_2)$ is a compact
length space.
\end{proposition}

Hereafter we assume that $(X,d)$ is a compact length space.
By definition, a {\em Wasserstein geodesic} is a minimizing geodesic in the
length space $(P(X), W_2)$. (We will always parametrize minimizing
geodesics in length spaces to have constant speed.)
The length space $(P(X), W_2)$ has some interesting features;
even for simple $X$, there may be an uncountable number of
Wasserstein geodesics between two measures $\mu_0, \mu_1 \in P(X)$
\cite[Example 2.9]{LV}.

As mentioned above, there is a relation between
minimizing geodesics in $(P(X), W_2)$  and minimizing geodesics in $X$.
Let $\Gamma$ be the set of minimizing
geodesics $\gamma \: : \: [0, 1] \rightarrow X$. 
It is compact in the uniform topology. 
For any $t\in [0,1]$, the {\em evaluation map} 
$e_t \: : \: \Gamma \rightarrow X$ defined by 
\begin{equation} \label{eqn3}
e_t(\gamma) \: = \: \gamma(t) \end{equation} 
is continuous. Let $E \: : \: \Gamma \rightarrow X \times X$ be the
``endpoints'' map given by $E(\gamma) \: = \: (e_0(\gamma),e_1(\gamma))$. 
A {\em dynamical transference plan} consists of a transference plan $\pi$ 
and a Borel measure $\Pi$ on $\Gamma$ such that $E_* \Pi \: = \: \pi$;
it is said to be optimal if $\pi$ itself is. 
In words, the transference plan $\pi$ tells us how much mass goes from
a point $x_0$ to another point $x_1$, but does not tell us about the actual
path that the mass has to follow. Intuitively, mass should flow
along geodesics, but there may be several possible choices of geodesics
between two given points and the transport may be divided among
these geodesics; this is the information provided
by $\Pi$. 

If $\Pi$ is an optimal dynamical transference plan then for $t \in [0,1]$, 
we put
\begin{equation} \label{eqn4}
\mu_t \: = \: (e_t)_* \Pi.
\end{equation} 
The one-parameter family of measures $\{\mu_t\}_{t \in [0,1]}$ 
is called a {\em displacement interpolation}.
In words, $\mu_t$ is what has become of the mass of $\mu_0$ after
it has travelled from time~0 to time~$t$ according to the dynamical
transference plan $\Pi$.

\begin{proposition} \cite[Lemma 2.4 and Proposition 2.10]{LV} \label{prop2}
Any displacement interpolation is a Wasserstein geodesic.
Conversely, any Wasserstein geodesic arises as a displacement
interpolation from some 
optimal dynamical transference plan.
\end{proposition}

In the Riemannian case, if $\mu_0, \mu_1$ are absolutely
continuous with respect to $\dvol_M$, and $F(m) \: = \:
\exp_m (- \: \nabla_m \phi)$ is the Monge transport 
between them, then there is a unique
Wasserstein geodesic between $\mu_0$ and $\mu_1$
given by $\mu_t \: = \: (F_t)_* \mu_0$, where
$F_t(m) \: = \: \exp_m (- \: t \: \nabla_m \phi)$.
Here $\mu_t$ is also absolutely continuous
with respect to $\dvol_M$. On the other hand, if $\mu_0 = \delta_{m_0}$ and 
$\mu_1 = \delta_{m_1}$ then some Wasserstein geodesics from
$\mu_0$ to $\mu_1$ are of the form
$\mu_t \: = \: \delta_{c(t)}$, where $c$ is a 
minimizing geodesic from $m_0$ to $m_1$. In particular,
the Wasserstein geodesic need not be unique.

In a remarkable paper \cite{Otto (2001)}, motivated by PDE problems,
Otto constructed a {\em formal}
infinite-dimensional
Riemannian metric $g_{H^{-1}}$ on $P(\R^n)$. To describe $g_{H^{-1}}$,
for simplicity we work with a compact
Riemannian manifold $M$ instead of $\R^n$. Suppose
that $\mu \in P(M)$ can be written as $\mu \: = \: \rho \:
\dvol_M$, with $\rho$ a smooth positive function. We formally think
of a tangent vector $\delta \mu \in T_\mu P(M)$ as being a
variation of $\mu$, which we take to be
$(\delta \rho) \: \dvol_M$ with $\delta \rho \in C^\infty(M)$.
There is a $\Phi \in C^\infty(M)$, unique up to constants,
so that $\delta \rho \: = \: d^* (\rho d \Phi)$. Then by
definition,
\begin{equation} \label{eqn5}
g_{H^{-1}}(\delta \mu, \delta \mu) \: = \:
\int_M |d \Phi|^2 \: d\mu.
\end{equation}
One sees that in terms of $\delta \rho \in C^\infty(M)$, $g_{H^{-1}}$
corresponds to a weighted $H^{-1}$-inner product.

Otto showed that the
corresponding distance function on $P(M)$ is {\em formally} $W_2$, and that
the ``infinite-dimensional Riemannian manifold'' $(P(\R^n),
g_{H^{-1}})$ {\em formally} has nonnegative sectional curvature.
One can make rigorous sense of these statements in terms of Alexandrov
geometry.

\begin{proposition} 
\cite[Theorem A.8]{LV},\cite[Proposition 2.10(iv)]{Sturm2} \label{prop3}
$(P(M), W_2)$ has nonnegative Alexandrov curvature if and only
if $M$ has nonnegative sectional curvature.
\end{proposition}

\begin{proposition} \cite[Proposition A.33]{LV} \label{prop4}
If $M$ has nonnegative sectional curvature then for each
absolutely continuous measure $\mu \: = \: \rho \: \dvol_M \in P(M)$, the
tangent cone $T_\mu P(M)$ is an inner product space. If $\rho$ is smooth
and positive
then the inner product on $T_\mu P(M)$ equals $g_{H^{-1}}$.
\end{proposition}

An open question is whether there is any good
sense in which $(P(M), W_2)$, or a large part thereof,
carries an infinite-dimensional Riemannian structure.
The analogous question for finite-dimensional Alexandrov
spaces has been much studied.

\begin{remark} \label{rmk2}
In Sturm's work he
uses the following interesting metric $D$ on the set of compact 
metric-measure spaces \cite[Definition 3.2]{Sturm2}. 
Given ${\mathcal X}_1 = (X_1, d_1, \nu_1)$ 
and ${\mathcal X}_2 = (X_2, d_2, \nu_2)$, let $\widehat{d}$ denote a
metric on the disjoint union $X_1 \coprod X_2$ such that
$\widehat{d} \big|_{X_1 \times X_1} \: = \: d_1$ and 
$\widehat{d} \big|_{X_2 \times X_2} \: = \: d_2$. Then
\begin{equation} \label{eqn6}
D({\mathcal X}_1, {\mathcal X}_2 )^2 \: = \:
\inf_{\widehat{d}, q} \int_{X_1 \times X_2}
\widehat{d}(x_1, x_2)^2 \: dq(x_1, x_2), 
\end{equation}
where $q$ runs over probability measures on $X_1 \times X_2$ whose
pushforwards onto $X_1$ and $X_2$ are $\nu_1$ and $\nu_2$, respectively.
If one restricts to metric-measure spaces with an upper diameter bound
whose measures have full
support and satisfy a uniform doubling condition
(which will be the case with a lower Ricci curvature bound)
then the topology coming from $D$ coincides with the MGH topology
of Definition \ref{def2} \cite[Lemma 3.18]{Sturm2},\cite{TOT2}.
\end{remark}

\section{Motivation for displacement convexity} \label{sec2}

To say a bit more about the PDE motivation, we recall that
the heat equation $\frac{\partial f}{\partial t} \: = \: \nabla^2 f$
can be considered to be the formal gradient flow of the
Dirichlet energy $E(f) \: = \: \frac12 \: \int_M |df|^2 \: \dvol_M$ on
$L^2(M, \dvol_M)$. (Our conventions are that a function decreases
along the flowlines of 
its gradient flow, so on a finite-dimensional Riemannian
manifold $Y$ the gradient
flow of a function $F \in C^\infty(Y)$ is $\frac{dc}{dt} \: = \:
- \: \nabla F$.) 
Jordan-Kinderlehrer-Otto showed that the heat equation on measures
can also be formally written as a gradient flow
\cite{Jordan-Kinderlehrer-Otto (1998)}.  Namely, for
a smooth
probability measure $\mu \: = \: \rho \: \frac{\dvol_M}{\vol(M)}$,
let us put $H_\infty(\mu) \: = \: \int_M \rho \: \log \rho \:
\frac{\dvol_M}{\vol(M)}$. Then the heat equation
$\frac{\partial }{\partial t} 
\left( \rho \: \frac{\dvol_M}{\vol(M)} \right) \: = \:
\nabla^2 \rho \: \frac{\dvol_M}{\vol(M)}$ 
is formally the gradient flow of $H_\infty$
on $P(M)$, where $P(M)$ has Otto's formal Riemannian metric.
Identifying a.c. measures and measurable functions using $\frac{\dvol_M}{\vol(M)}$, this
gave a new way to realize the heat equation as a gradient flow.

Although this approach may not give much new
information about the heat equation,
it has more relevance if one considers other functions $H$ on $P(M)$,
whose gradient flows can give rise to interesting nonlinear
PDE's such as the porous medium equation. Again formally,
if one has positive lower bounds on the Hessian of $H$ then one can 
draw conclusions about uniqueness of critical points and 
rates of convergence of the gradient flow to the critical point,
which one can then hope to make rigorous.
This reasoning motivated McCann's notion of {\em displacement convexity},
i.e. convexity of a function $H$ along Wasserstein geodesics 
\cite{McCann (1997)}.
(We recall that on a smooth manifold, a smooth function
has a nonnegative Hessian if and only
if it is convex when restricted to each geodesic.)

In a related direction, Otto and Villani \cite{Otto-Villani (2000)} saw that
convexity properties on $P(M)$ could be used to give heuristic arguments for
functional inequalities on $M$, such as the log Sobolev inequality.  They
could then give rigorous proofs based on these heuristic arguments.
Given a smooth background probability measure
$\nu \: = \: e^{- \Psi} \: \dvol_M$ and an absolutely 
continuous probability measure
$\mu \: = \: \rho \: \nu$, let us now put
$H_\infty(\mu) \: = \: \int_M \rho \: (\log \rho) \: d\nu$. 
As part of their work,
Otto and Villani computed the formal Hessian of the
function $H_\infty$ on $P(M)$
and found that it is bounded below
by $K g_{H^{-1}}$ provided that the Bakry-\'Emery tensor
$\Ric_\infty \: = \: \Ric \: + \: \Hess(\Psi)$ satisfies
$\Ric_\infty \: \ge \: K g$ on $M$. 
This was perhaps the first indication that Ricci curvature
is related to convexity properties on Wasserstein space.

Around the same time,
Cordero-Erausquin-McCann-Schmuckenschl\"ager \cite{CEMS01}
gave a rigorous proof of the convexity of certain functions
on $P(M)$ when $M$ has dimension $n$ and nonnegative Ricci
curvature.  Suppose that $A \: : \: [0, \infty) \rightarrow
\R$ is a continuous convex function
with $A(0) = 0$ such that $\lambda \rightarrow 
\lambda^n A(\lambda^{-n})$ is a convex  function on $\R^+$.
If $\mu \: = \: \rho \: \frac{\dvol_M}{\vol(M)}$ is  an
absolutely continuous probability measure then 
put $H_A(\mu) \: = \: \int_M A(\rho) \: \frac{\dvol_M}{\vol(M)}$. The
statement is that if 
$\mu_0, \mu_1 \in P(M)$ are absolutely continuous, and
$\{\mu_t\}_{t \in [0,1]}$ is the (unique) Wasserstein 
geodesic between them, then $H_A(\mu_t)$ is convex in $t$,
again under the assumption of nonnegative Ricci curvature.

Finally, von Renesse and Sturm \cite{Sturm-von Renesse} extended the work of
Cordero-Erausquin-McCann-Schmuckenschl\"ager to show that
the function $H_\infty$, defined by
$H_\infty \left(\rho \: \frac{\dvol_M}{\vol(M)} \right) \: = \: \int_M \rho \:
\log \rho \: \frac{\dvol_M}{\vol(M)}$, 
is $K$-convex along Wasserstein geodesics 
between  absolutely-continuous measures if and only if $\Ric \: \ge \: Kg$.
(The relation with the Otto-Villani result is that
$\Psi$ is taken to be constant, so $\nu \: = \: \frac{\dvol_M}{\vol(M)}$.)
The ``if'' implication is along the lines of the
Cordero-Erausquin-McCann-Schmuckenschl\"ager result and the
``only if'' implication involves some local arguments.

Although these results indicate a formal relation between 
Ricci curvature and
displacement convexity, one can ask for a more intuitive
understanding. Here is one example.
\begin{example}
Consider the functional
$H_\infty \left(
\rho \: \frac{\dvol_M}{\vol(M)} \right) \: = \: \int_M \rho \:
\log \rho \: \frac{\dvol_M}{\vol(M)}$. It is minimized, among
absolutely continuous probability measures on $M$, when
$\rho \: = \: 1$, i.e. when the measure 
$\mu \: = \: \rho \: \frac{\dvol_M}{\vol(M)}$ is the
uniform measure $\frac{\dvol_M}{\vol(M)}$. In this sense, $H_\infty$ measures
the {\em nonuniformity} of $\mu$ with respect to
$\frac{\dvol_M}{\vol(M)}$. Now take $M = S^2$. Let $\mu_0$ and 
$\mu_1$ be two small congruent rotationally symmetric
blobs, centered at the
north and south poles respectively.  Clearly
$U_\infty(\mu_0) \: = \: U_\infty(\mu_1)$. Consider the
Wasserstein geodesic from $\mu_0$ to $\mu_1$. It takes
the blob $\mu_0$ and pushes it down in a certain way 
along the lattitudes until it becomes $\mu_1$. At an
intermediate time, say around $t = \frac12$, the blob
has spread out to form a ring. When it spreads, it becomes
{\em more} uniform with respect to $\frac{\dvol_M}{\vol(M)}$. Thus the
{\em nonuniformity} at an intermediate time is {\em at most}
that at times $t =0$ or $t = 1$. This can be seen as a
consequence of the convexity of $H_\infty(\mu_t)$ in $t$,
i.e. for $t \in [0,1]$ we have
$H_\infty(\mu_t) \: \le H_\infty(\mu_0) \: = \:
H_\infty(\mu_1)$. In this way the displacement convexity
of $H_\infty$ can be seen as an averaged form of the focusing
property of positive curvature.  Of course this
example does not indicate why the relevant curvature is
Ricci curvature, as opposed to some other curvature, but
perhaps gives some indication of why curvature is related
to displacement convexity.
\end{example}

\section{Entropy functions and displacement convexity} \label{sec3}

In this section we give the definition of nonnegative
$N$-Ricci curvature. We then outline the proof that it is
preserved under measured Gromov-Hausdorff limits.  In the
next section we relate the definition to the classical notion
of Ricci curvature, in the case of a smooth metric-measure space.

\subsection{Definitions}

We first define the relevant ``entropy'' functionals.

Let $X$ be a compact Hausdorff space.
Let $U \: : \: [0, \infty) \rightarrow \R$ be a continuous
convex function with $U(0) \: = \: 0$.
Given a reference probability measure $\nu \in P(X)$,
define the entropy function
$U_\nu \: : \: P(X) \rightarrow \R \cup \{\infty\}$ by
\begeq \label{eqn6.5}
U_\nu(\mu) = 
\int_X U(\rho(x))\,d\nu(x) + U'(\infty)\,\mu_s(X), 
\endeq
where 
\begin{equation} \label{eqn7}
\mu = \rho \nu + \mu_s 
\end{equation}
is the Lebesgue decomposition of $\mu$ with respect to $\nu$ 
into an absolutely continuous part $\rho\nu$ and a singular part $\mu_s$, and
\begin{equation} \label{eqn8}
U^\prime(\infty) \: = \: \lim_{r \rightarrow \infty} \frac{U(r)}{r}.
\end{equation}

\begin{example}
Given $N \in (1, \infty]$, take the function 
$U_N$ on $[0, \infty)$ to be
\begin{equation} \label{eqn9}
U_N(r) \: = \: 
\begin{cases}
Nr (1 - r^{-1/N})& \text{ if $1 < N < \infty$}, \\
r \log{r} & \text{ if $N = \infty$}.
\end{cases}
\end{equation}
Let
$H_{N,\nu} \: : \: P(X) \rightarrow \R \cup \{\infty \}$ be
the corresponding entropy function. If
$N \in (1, \infty)$ then
\begin{equation} \label{eqn10}
H_{N,\nu} = N - N \int_X \rho^{1-\frac{1}{N}} \,d\nu, 
\end{equation} 
while if $N = \infty$ then
\begin{equation} \label{eqn11}
H_{\infty,\nu}(\mu) = \int_X \rho\log\rho \,d\nu
\end{equation}
if $\mu$ is absolutely continuous with
respect to $\nu$ and $H_{\infty,\nu}(\mu) = \infty$ otherwise.
\end{example}

One can show that as a function of $\mu \in P(X)$, $U_\nu(\mu)$ is
minimized when $\mu \: = \: \nu$. It would be better to
call $U_\nu$ a ``negative entropy'', but we will be sloppy.
Here are the technical properties of $U_\nu$ that we need.

\begin{proposition} \cite{Liese-Vajda (1987)},\cite[Theorem B.33]{LV} 
\label{prop5}
(i) $U_\nu(\mu)$ is a lower semicontinuous function of $(\mu,\nu)
\in P(X) \times P(X)$.
That is, if $\{\mu_k\}_{k=1}^\infty$ and 
$\{\nu_k\}_{k=1}^\infty$ are sequences in $P(X)$
with $\lim_{k \rightarrow \infty} \mu_k \: = \: \mu$ and
$\lim_{k \rightarrow \infty} \nu_k \: = \: \nu$ in the weak-$*$ topology then
\begin{equation} \label{eqn12}
U_\nu(\mu) \leq \liminf_{k\to\infty} U_{\nu_k}(\mu_k).
\end{equation}

\quad (ii) $U_\nu(\mu)$ is nonincreasing under pushforward.
That is, if $Y$ is a compact Hausdorff space and
$f \: : \: X \rightarrow Y$ is a Borel map then
\begin{equation} \label{eqn13}
U_{f_* \nu} (f_*\mu) \leq U_\nu(\mu).
\end{equation}
\end{proposition}

In fact, the $U^\prime(\infty) \: \mu_s(X)$ term in (\ref{eqn6.5}) is 
dictated by the fact that we want
$U_\nu$ to be lower semicontinuous on $P(X)$.

We now pass to the setting of a compact measured length space $(X, d, \nu)$.
The definition of nonnegative $N$-Ricci curvature will be
in terms of the convexity of certain
entropy functions on $P(X)$, where the entropy is relative to
the background measure $\nu$. By ``convexity'' we mean
convexity along Wasserstein geodesics, i.e. displacement
convexity. We first describe the
relevant class of entropy functions.
 
If $N\in [1,\infty)$ then 
we define $\DC_N$ to be the set of such functions $U$
so that the function
\begin{equation} \label{eqn14}
\psi(\lambda) = \lambda^N \:U(\lambda^{-N}) 
\end{equation}
is convex on $(0, \infty)$.
We further define $\DC_\infty$ 
to be the set of such functions $U$ so that the function
\begin{equation} \label{eqn15}
\psi(\lambda) = e^\lambda \: U(e^{-\lambda})
\end{equation} 
is convex on $(- \infty, \infty)$.
A relevant example of an element of $\DC_N$ is given by
the function $U_N$ of (\ref{eqn9}).

\begin{definition} \cite[Definition 5.12]{LV} \label{def3}
Given $N \in [1, \infty]$, we say that a compact measured length space
$(X,d,\nu)$ has nonnegative $N$-Ricci curvature if for all
$\mu_0, \mu_1 \in P(X)$ with
$\supp(\mu_0) \subset \supp(\nu)$ and $\supp(\mu_1) \subset \supp(\nu)$, 
there is {\em some} Wasserstein geodesic
$\{\mu_t\}_{t \in [0,1]}$ from $\mu_0$ to $\mu_1$ so that
for all $U \in \DC_N$ and all $t \in [0,1]$,
\begin{equation} \label{ineq1}
U_\nu(\mu_t) \: \le \: t \: U_\nu(\mu_1) \: + \: (1-t) \: 
U_\nu(\mu_0).
\end{equation}
\end{definition}

We make some remarks about the definition.
\begin{remark} \label{rmk3}
A similar definition in the case $N = \infty$, but in terms
of $U = U_\infty$ instead of $U \in \DC_\infty$, was used in
\cite[Definition 4.5]{Sturm2}; see also Remark \ref{rmk9}.
\end{remark}
\begin{remark} \label{rmk4}
It is not hard to show that if $(X,d,\nu)$ has nonnegative
$N$-Ricci curvature and $N^\prime \ge N$ then
$(X,d,\nu)$ has nonnegative $N^\prime$-Ricci curvature. 
\end{remark}
\begin{remark} \label{rmk5}
Note that for $t \in (0,1)$,
the intermediate measures $\mu_t$ are not required to have support
in $\supp(\nu)$. If $(X,d,\nu)$ has nonnegative $N$-Ricci curvature
then $\supp(\nu)$ is a convex subset of $X$ and
$(\supp(\nu), d \big|_{\supp(\nu)}, \nu)$ has nonnegative $N$-Ricci
curvature \cite[Theorem 5.53]{LV}.
(We recall that a subset $A \subset X$ is {\em convex} if for 
any $x_0, x_1 \in A$ there is a minimizing geodesic from $x_0$ to $x_1$ 
that lies entirely in $A$. It is {\em totally convex} if for any
$x_0, x_1 \in A$, any minimizing geodesic in $X$ from $x_0$ to
$x_1$ lies in $A$.)
So we don't lose much by assuming that
$\supp(\nu) \: = \: X$.
\end{remark}
\begin{remark} \label{rmk6}
There is supposed to be a single Wasserstein geodesic
$\{\mu_t\}_{t \in [0,1]}$ from $\mu_0$ to $\mu_1$ so that
(\ref{ineq1}) holds along $\{\mu_t\}_{t \in [0,1]}$ for all
$U \in \DC_N$ simultaneously. However,
(\ref{ineq1}) is only assumed to
hold along {\em some} Wasserstein geodesic from $\mu_0$ to $\mu_1$, 
and not necessarily along all such Wasserstein geodesics.  This is what we call
{\em weak displacement convexity}. It may be more conventional
to define convexity on a length space in terms of convexity along
{\em all} geodesics. However, the definition
with weak displacement convexity turns out to work better under MGH limits, 
and has most of the same implications as if we required
convexity along all Wasserstein geodesics from $\mu_0$ to $\mu_1$.
\end{remark}
\begin{remark} \label{rmk7}
Instead of requiring that (\ref{ineq1}) holds for all $U \in \DC_N$,
it would be consistent to make a definition in which it is only
required to hold for the function $U = U_N$ of (\ref{eqn9}).
For technical reasons, we
prefer to require that (\ref{ineq1}) holds for all $U \in \DC_N$;
see Remark \ref{rmk13}. Also,
the class $\DC_N$ is the natural class of
functions for which the proof of Theorem \ref{thm2} works.
\end{remark}

\subsection{MGH invariance}

The next result says that Definition \ref{def3} satisfies
Condition 1. of Wishlist \ref{wishlist}. It shows that for each $N$,
there is a self-contained world of measured length spaces with
nonnegative $N$-Ricci curvature.

\begin{theorem} 
\cite[Theorem 5.19]{LV},\cite[Theorem 4.20]{Sturm2},
\cite[Theorem 3.1]{Sturm3} 
\label{thm1}
Let $\{(X_i, d_i, \nu_i)\}_{i=1}^\infty$ be a sequence of
compact measured length spaces with $\lim_{i \rightarrow \infty} 
(X_i,d_i,\nu_i) \: = \: (X,d,\nu)$ in the
measured Gromov-Hausdorff topology.
For any $N \in [1, \infty]$, if 
each $(X_i, d_i, \nu_i)$ has nonnegative $N$-Ricci curvature then
$(X, d, \nu)$ has nonnegative $N$-Ricci curvature.
\end{theorem}
\begin{proof}
We give an outline of the proof.  For simplicity, we just
consider a single $U \in \DC_N$; the same argument will allow one
to handle all $U \in \DC_N$ simultaneously.

Suppose first that $\mu_0$ and $\mu_1$ are absolutely continuous
with respect to $\nu$, with continuous densities $\rho_0, \rho_1
\in C(X)$. 
Let $\{f_i\}_{i=1}^\infty$ be a sequence of $\epsilon_i$-approximations
as in Definition \ref{def2}.
We first approximately-lift the measures $\mu_0$ and $\mu_1$ to $X_i$. 
That is, we use $f_i$ to pullback the densities to $X_i$, then
multiply by $\nu_i$ and
then normalize to get probability measures.
More precisely, we put
$\mu_{i,0} \: = \: \frac{f_i^* \rho_0 \: \nu_i}{
\int_{X_i} f_i^* \rho_0 \: d\nu_i} \in P(X_i)$ and 
$\mu_{i,1} \: = \: \frac{f_i^* \rho_1 \: \nu_i}{
\int_{X_i} f_i^* \rho_0 \: d\nu_i} \in P(X_i)$. One shows that
$\lim_{i \rightarrow \infty} (f_i)_* \mu_{i,0} \: = \: \mu_{0}$ and
$\lim_{i \rightarrow \infty} (f_i)_* \mu_{i,1} \: = \: \mu_{1}$ in
the weak-$*$ topology on $P(X)$. In addition, one shows that
\begin{equation} \label{eqn16}
\lim_{i \rightarrow \infty} U_{\nu_i}(\mu_{i,0}) \: = \: U_\nu(\mu_0)
\end{equation} 
and
\begin{equation} \label{eqn17}
\lim_{i \rightarrow \infty} U_{\nu_i}(\mu_{i,1}) \: = \: U_\nu(\mu_1).
\end{equation}

Up on $X_i$, we are OK in the sense that
by hypothesis, there is a Wasserstein geodesic
$\{\mu_{i,t}\}_{t \in [0,1]}$ from 
$\mu_{i,0}$ to $\mu_{i,1}$ in $P(X_i)$ so that for all $t \in [0,1]$, 
\begin{equation} \label{eqn18}
U_{\nu_i}(\mu_{i,t}) \: \le \: t \: U_{\nu_i}(\mu_{i,1}) \: + \: (1-t) \: 
U_{\nu_i}(\mu_{i,0}).
\end{equation}
We now want to take a convergent subsequence
of these Wasserstein geodesics in an appropriate sense 
to get a Wasserstein geodesic in $P(X)$. This can be done
using Proposition \ref{prop1} and an Arzela-Ascoli-type result.
The conclusion is that after
passing to a subsequence of the $i$'s, there is a Wasserstein geodesic  
$\{\mu_{t}\}_{t \in [0,1]}$ from $\mu_{0}$ to $\mu_{1}$ in $P(X)$
so that for each $t \in [0,1]$, we have $\lim_{i \rightarrow \infty}
(f_i)_* \mu_{i,t} \: = \: \mu_t$. 

Finally, we want to see what (\ref{eqn18}) becomes as $i \rightarrow \infty$.
At the endpoints we have good limits from (\ref{eqn16}) and (\ref{eqn17}), so
this handles the right-hand-side of (\ref{eqn18}) as $i \rightarrow \infty$.
We do not have such a good limit for the left-hand-side. However, this
is where the lower semicontinuity comes in.
Applying parts (i) and (ii) of Proposition \ref{prop5}, we do know that
\begin{equation} \label{eqn19}
U_{\nu}(\mu_{t}) \: \le \: \liminf_{i \rightarrow \infty} 
U_{(f_i)_* \nu_i}((f_i)_* \mu_{i,t})
 \: \le \: \liminf_{i \rightarrow \infty} 
U_{\nu_i}(\mu_{i,t}).
\end{equation}
This is enough to give the desired inequality (\ref{ineq1}) along
the Wasserstein geodesic $\{\mu_{t}\}_{t \in [0,1]}$.

This handles the case when $\mu_0$ and $\mu_1$ have continuous densities.
For general $\mu_0, \mu_1 \in P(X)$, using mollifiers we can construct
sequences $\{ \mu_{j,0} \}_{j=1}^\infty$
and $\{ \mu_{j,1} \}_{j=1}^\infty$ 
of absolutely continuous measures with continuous densities so that
$\lim_{j \rightarrow \infty} \mu_{j,0} \: = \: \mu_0$ and
$\lim_{j \rightarrow \infty} \mu_{j,1} \: = \: \mu_1$ in the
weak-$*$ topology.  In addition, one can do the mollifying in such a way
that
$\lim_{j \rightarrow \infty} U_\nu(\mu_{j,0}) \: = \: U_\nu(\mu_0)$ and
$\lim_{j \rightarrow \infty} U_\nu(\mu_{j,1}) \: = \: U_\nu(\mu_1)$.
From what has already been shown, for each $j$ there is
a Wasserstein geodesic $\{\mu_{j,t}\}_{t \in [0,1]}$ in
$P(X)$ from $\mu_{j,0}$ to $\mu_{j,1}$ so that
for all $t \in [0,1]$, 
\begin{equation} \label{eqn20}
U_{\nu}(\mu_{j,t}) \: \le \: t \: U_{\nu}(\mu_{j,1}) \: + \: (1-t) \: 
U_{\nu}(\mu_{j,0}).
\end{equation}
After passing to a subsequence, we can assume that the
Wasserstein geodesics $\{\mu_{j,t}\}_{t \in [0,1]}$ converge
uniformly as $j \rightarrow \infty$ to a Wasserstein geodesic 
$\{\mu_{t}\}_{t \in [0,1]}$ from $\mu_0$ to $\mu_1$. From the
lower semicontinuity of $U_\nu$, we have
$U_{\nu}(\mu_{t}) \: \le \: \liminf_{j \rightarrow \infty} 
U_{\nu}(\mu_{j,t})$. Equation (\ref{ineq1}) follows.
\end{proof}

\subsection{Basic properties}

We now give some basic properties of measured length spaces 
$(X,d,\nu)$ with
nonnegative $N$-Ricci curvature.

\begin{proposition} \cite[Proposition 5.20]{LV},
\cite[Theorem 2.3]{Sturm3} \label{prop6}
For $N \in (1, \infty]$, if $(X,d,\nu)$  has
nonnegative $N$-Ricci curvature then
the measure $\nu$ is either a delta function or is nonatomic.
The support of $\nu$ is a convex subset of $X$.
\end{proposition}

The next result is an analog of the Bishop-Gromov theorem.

\begin{proposition} \cite[Proposition 5.27]{LV},
\cite[Theorem 2.3]{Sturm3} \label{prop7}
Suppose that $(X,d,\nu)$ has nonnegative $N$-Ricci curvature, with
$N \in [1,  \infty)$. Then for all 
$x\in \supp(\nu)$ and all $0<r_1\leq r_2$,
\begin{equation} \label{eqn21}
\nu(B_{r_2}(x)) \leq \left ( \frac{r_2}{r_1} \right )^N \nu(B_{r_1}(x)).
\end{equation}
\end{proposition}
\begin{proof}
We give an outline of the proof.
There is a Wasserstein
geodesic $\{ \mu_t\}_{t \in [0,1]}$ between 
$\mu_0 \: = \: \delta_{x}$ and the restricted measure
$\mu_1 \: = \: \frac{1_{B_{r_2}(x)}}{\nu(B_{r_2}(x))} \: \nu$,
along which (\ref{ineq1}) holds. Such a Wasserstein geodesic comes
from a fan of geodesics (the support of $\Pi$)
that go from $x$ to points in $B_{r_2}(x)$.
The actual transport, going backwards from $t=1$ to $t=0$, amounts to
sliding the mass of $\mu_1$ along these geodesics towards
$x$. In particular, the support of $\mu_t$ is contained in
$B_{tr_2}(x)$. Applying (\ref{ineq1}) with $U = U_N$ and
$t \: = \: \frac{r_1}{r_2}$, along with Holder's inequality,
gives the desired result.
\end{proof}

We give a technical result which will be used in
deriving functional inequalities.

\begin{proposition} \cite[Theorem 5.52]{LV} \label{prop8}
Suppose that $(X, d, \nu)$ has nonnegative $N$-Ricci curvature.
If $\mu_0$ and $\mu_1$ are absolutely continuous with respect to
$\nu$ then the measures in the Wasserstein geodesic 
$\{\mu_t\}_{t \in [0,1]}$ of Definition \ref{def3} are all
absolutely continuous with respect to $\nu$.
\end{proposition}

Finally, we mention that for {\em nonbranching} measured
length spaces, there is a local-to-global principle which 
says that having nonnegative $N$-Ricci curvature in a local
sense implies nonnegative $N$-Ricci curvature in a global
sense \cite[Theorem 4.17]{Sturm2},\cite{TOT2}.
We do not know if this holds in the branching case.

\section{Smooth metric-measure spaces} \label{sec4}

We now address Condition 2. of Wishlist \ref{wishlist}.
We want to know what our abstract definition of
``nonnegative $N$-Ricci curvature'' boils down to in the
classical Riemannian case.  To be a bit more general, we
allow Riemannian manifolds with weights. Let
us say that a {\em smooth} measured length space consists
of a smooth $n$-dimensional Riemannian manifold $M$ along with a smooth
probability measure $\nu \:  = \: e^{- \: \Psi} \: \dvol_M$.
We write $(M,g,\nu)$ for the corresponding measured length
space. We are taking $M$ to be compact.

Let us discuss possible Ricci tensors for smooth measured
length spaces.
If $\Psi$ is constant, i.e. if $\nu \: = \: \frac{\dvol_M}{\vol(M)}$, then
the right notion of a Ricci tensor for $M$ is clearly just the usual
$\Ric$. 

For general $\Psi$, a modified Ricci tensor 
\begin{equation} \label{eqn22}
\Ric_\infty \: = \: \Ric \: + \: \Hess(\Psi)
\end{equation} 
was introduced
by Bakry and \'Emery \cite{Bakry-Emery (1985)}. 
(Note that the standard $\R^n$ with the Gaussian measure
$(2\pi)^{-\frac{n}{2}} \: e^{- \frac{|x|^2}{2}} \: d^nx$ has a
constant Bakry-\'Emery tensor given by
$(\Ric_\infty)_{ij} \: = \: \delta_{ij}$.) 
Their motivation came from a desire
to generalize the Lichnerowicz inequality for the lower positive
eigenvalue $\lambda_1(\triangle)$ 
of the Laplacian.  We recall the Lichnerowicz result
that if an $n$-dimensional Riemannian manifold has 
$\Ric \: \ge \: K \: g$ with $K > 0$  then
$\lambda_1(\triangle) \: \ge \: \frac{n}{n-1} \: K$
\cite{Lichnerowicz}. 

In the case of a 
Riemannian manifold with a smooth probability measure $\nu \: = \: 
e^{- \: \Psi} \: \dvol_M$, there is a natural self-adjoint Laplacian
$\widetilde{\triangle}$ acting on the weighted $L^2$-space
$L^2(M, e^{- \: \Psi} \: \dvol_M)$, given by
\begin{equation} \label{eqn23}
\int_M f_1 (\widetilde{\triangle} f_2) \: e^{- \: \Psi} \: \dvol_M \: = \:
\int_M  \langle \nabla f_1, \nabla f_2 \rangle \: e^{- \: \Psi} \: \dvol_M
\end{equation}
for $f_1, f_2 \in C^\infty(M)$. Here $\langle \nabla f_1, \nabla f_2 \rangle$
is the usual local inner product
computed using the Riemannian metric $g$.
Bakry and \'Emery showed that if $\Ric_\infty \: \ge \: K g$ then
$\lambda_1 (\widetilde{\triangle}) \: \ge \: K$. Although this
statement is missing the $\frac{n}{n-1}$ factor of the
Lichnerowicz inequality, it holds independently of $n$ and so
can be considered to be a version of the Lichnerowicz inequality
where one allows weights and takes $n \rightarrow \infty$. We
refer to \cite{Toulouse (2000)} 
for more information on the Bakry-\'Emery tensor
$\Ric_\infty$,
including its relationship to log Sobolev inequalities. 
Some geometric properties of $\Ric_\infty$ were studied in
\cite{Lott (2003)}. More
recently, the Bakry-\'Emery tensor has appeared as the 
right-hand-side of Perelman's modified Ricci flow equation
\cite{Perelman}.

We have seen that $\Ric_\infty$ is a sort of Ricci tensor
for the smooth measured length space $(M,g,\nu)$ when we
consider $(M,g,\nu)$ to have ``effective dimension'' infinity.
There is a similar tensor for other effective dimensions.
Namely, if $N \in (n, \infty)$ then we put
\begin{equation} \label{eqn24}
\Ric_N \: = \: \Ric \: + \: \Hess(\Psi) \: - \: 
\frac{1}{N-n} \: d\Psi \otimes d\Psi,
\end{equation}
where $\dim(M) = n$. The intuition is that $(M,g,\nu)$ has
conventional dimension $n$ but is pretending to have
dimension $N$, and $\Ric_N$ is its effective Ricci tensor
under this pretence. There is now a sharp analog of the
Lichnerowicz inequality : if $\Ric_N \: \ge \: K g$ with
$K > 0$ then $\lambda_1(\widetilde{\triangle}) \: \ge \:
\frac{N}{N-1} \: K$ \cite{Bakry}. Geometric properties of
$\Ric_N$ were studied in \cite{Lott (2003)} and 
\cite{Qian (1997)}.

Finally, if $N < n$, or if $N = n$ and $\Psi$ is not
locally constant, then we take the effective Ricci
tensor $\Ric_N$ to be $- \infty$. To summarize,
\begin{definition} \label{def4}
For $N \in [1, \infty]$, define the $N$-Ricci tensor $\Ric_N$ of $(M, g, \nu)$
by
\begin{equation} \label{eqn25}
\Ric_N \: = 
\begin{cases}
\Ric \: + \: \Hess(\Psi)  & \text{ if $N = \infty$}, \\
\Ric \: + \: \Hess(\Psi) \: - \: \frac{1}{N-n} \: d \Psi \otimes
d \Psi & \text{ if $n \: < \: N \: < \: \infty$}, \\
\Ric \: + \: \Hess(\Psi) \: - \: \infty \: (d \Psi \otimes
d \Psi) & \text{ if $N = n$}, \\
- \infty & \text{ if $N < n$,}
\end{cases}
\end{equation}
where by convention $\infty \cdot 0 \: = \: 0$.
\end{definition}

We can now state what the abstract notion of 
nonnegative $N$-Ricci curvature boils down to in the
smooth case.

\begin{theorem} 
\cite[Theorems 7.3 and 7.42]{LV},\cite[Theorem 4.9]{Sturm2},
\cite[Theorem 1.7]{Sturm3} \label{thm2}
Given $N \in [1, \infty]$, the measured length space 
$(M, g, \nu)$ has nonnegative $N$-Ricci curvature 
in the sense of Definition \ref{def3} if and
only if $\Ric_N \: \ge \: 0$.
\end{theorem}

The proof of Theorem \ref{thm2} uses the explicit description of
optimal transport on Riemannian manifolds.

In the special case when $\Psi$ is constant, and so
$\nu \: = \: \frac{\dvol_M}{\vol(M)}$,
Theorem \ref{thm2} shows that we recover the usual notion
of nonnegative Ricci curvature from our length space definition
as soon as $N\geq n$.

\subsection{Ricci limit spaces}

We give an application of Theorems \ref{thm1} and \ref{thm2} 
to Ricci limit spaces. 
From Gromov precompactness, given $N \in \Z^+$ and $D > 0$, the
Riemannian manifolds with
nonnegative Ricci curvature, dimension at most $N$ and
diameter at most $D$ form a precompact subset of the set of
measured length spaces, with respect to the MGH topology. 
The problem is to characterize the limit points.  In general
the limit points can be very singular, so this is a hard problem.
However, let us ask a simpler question : what are the 
limit points that happen to be {\em smooth} measured length
spaces?  That is, we are trying to characterize the smooth
limit points.

\begin{corollary} \cite[Corollary 7.45]{LV} \label{maincor}
If $(B, g_B, e^{-  \Psi} \dvol_B)$ is a measured 
Gromov-Hausdorff limit of Riemannian manifolds with nonnegative
Ricci curvature and dimension at most $N$ then $\Ric_N(B) \: \ge \:0$. 
(Here $B$ has dimension $n$, which is less than or equal to $N$.)
\end{corollary}
\begin{proof}
Suppose that $\{(M_i, g_i)\}_{i=1}^\infty$ is a sequence of
Riemannian manifolds with nonnegative Ricci curvature and
dimension at most $N$, with
$\lim_{i \rightarrow \infty} \left( M_i, g_i, 
\frac{\dvol_{M_i}}{\vol(M_i)} \right)
\: = \: (B, g_B, e^{-  \Psi} \dvol_B)$. From Theorem \ref{thm2}, the measured
length space $\left( M_i, g_i, 
\frac{\dvol_{M_i}}{\vol(M_i)} \right)$ has nonnegative $N$-Ricci
curvature.  From Theorem \ref{thm1}, $(B, g_B, e^{-  \Psi} \dvol_B)$ has
nonnegative $N$-Ricci curvature.  From Theorem \ref{thm2} again,
$\Ric_N(B) \: \ge \: 0$.
\end{proof}

There is a partial converse to Corollary \ref{maincor}.
\begin{proposition} \cite[Corollary 7.45]{LV} \label{prop9}
(i) Suppose that $N$ is an integer. If $(B, g_B, e^{-  \Psi} \dvol_B)$ 
has $\Ric_N(B) \: \ge \: 0$ with $N \: \ge \: \dim(B) \: + \: 2$ then 
$(B,  g_B, e^{-  \Psi} \dvol_B)$ is a measured 
Gromov-Hausdorff limit of Riemannian manifolds with nonnegative
Ricci curvature and dimension $N$. \\
(ii) Suppose that $N = \infty$. If $(B, g_B, e^{-  \Psi} \dvol_B)$ has 
$\Ric_\infty(B) \: \ge \: 0$
then $(B, g_B, e^{-  \Psi} \dvol_B)$ is a measured 
Gromov-Hausdorff limit of Riemannian manifolds $M_i$ with 
$\Ric(M_i) \: \ge \: - \: \frac{1}{i} \: g_{M_i}$.
\end{proposition}
\begin{proof}
Let us consider part (i). The proof uses the warped
product construction of \cite{Lott (2003)}. Let $g_{S^{N-\dim(B)}}$ be the
standard metric on the sphere $S^{N-\dim(B)}$. Let $M_i$ be
$B \times S^{N-\dim(B)}$ with the warped product metric
$g_i \: = \: g_B \: + \: i^{-2} e^{- \: \frac{\Psi}{N-\dim(B)}} \:
g_{S^{N-\dim(B)}}$. The metric $g_i$ is constructed so that
if $p \: : \: B \times S^{N-\dim(B)} \rightarrow B$ is projection
onto the first factor then $p_* \dvol_{M_i}$ is a constant
times $e^{-  \Psi} \dvol_B$. 
In terms of the fibration $p$, the Ricci tensor of 
$M_i$ splits into
horizontal and vertical components, with the horizontal
component being exactly $\Ric_N$.
As $i$ increases, the fibers shrink and the vertical
Ricci curvature of $M_i$ becomes dominated by the
Ricci curvature of the small fiber $S^{N-\dim(B)}$, which is positive
as we are assuming that $N-\dim(B) \ge 2$. Then for large $i$,
$(M_i, g_i)$ has nonnegative Ricci curvature.  Taking
$f_i = p$, we see that $\lim_{i \rightarrow \infty} \left( M_i, g_i, 
\frac{\dvol_{M_i}}{\vol(M_i)} \right)
\: = \: (B, g_B, e^{-  \Psi} \dvol_B)$.

The proof of (ii) is similar, except that we also allow
the dimensions of the fibers to go to infinity.
\end{proof}

Examples of singular spaces with nonnegative $N$-Ricci curvature
come from group actions.
Suppose that a compact Lie group $G$ acts isometrically on a 
$N$-dimensional Riemannian manifold $M$ that has
nonnegative Ricci curvature. Put $X = M/G$, let
$p \: : \: M \rightarrow X$ be the quotient map, let
$d$ be the quotient metric and
put $\nu \: = \: p_* \left( 
\frac{\dvol_M}{\vol(M)} \right)$. Then
$(X, d, \nu)$ has nonnegative $N$-Ricci curvature
\cite[Corollary 7.51]{LV}.

Finally, we recall the theorem of O'Neill that sectional
curvature is nondecreasing under pushforward by a Riemannian
submersion.  There is a Ricci analog of the O'Neill
theorem, expressed in terms of the modified Ricci tensor $\Ric_N$
\cite{Lott (2003)}.
The proof of this in 
\cite{Lott (2003)} was by explicit tensor calculations.
Using optimal transport, one can give a ``synthetic'' proof
of this Ricci O'Neill theorem 
\cite[Corollary 7.52]{LV}.  (This is what first convinced
the author that optimal transport is the right approach.)

\begin{remark} \label{rmk8}
We return to the question of whether one can give a good definition of
``nonnegative $N$-Ricci curvature'' by just taking the conclusion
of the Bishop-Gromov theorem and turning it into a definition.
To be a bit more reasonable, we consider taking an angular Bishop-Gromov 
inequality as the definition. Such an
inequality, with parameter $n$, does indeed characterize when
an $n$-dimensional Riemannian manifold has nonnegative Ricci curvature.
Namely, from comparison geometry, nonnegative Ricci curvature implies an
angular Bishop-Gromov inequality. To go the other way, suppose that
the angular Bishop-Gromov inequality holds.
We use polar coordinates around a point $m \in M$ and
recall that the volume of a infinitesimally small angular sector 
centered in the direction of a unit vector 
$v \in T_mM$, and going up to radius $r$, has the
Taylor expansion 
\begin{equation} \label{eqn26}
V(v, r) \: = \: \const \: r^n \: \left( 1 \: - \: 
\frac{n}{6(n+2)} \: \Ric(v,v) \: r^2 \: + \: \ldots \right).
\end{equation}
If $r^{-n} \: V(v,r)$ is to be nonincreasing in $r$ then we must have
$\Ric(v,v) \: \ge \: 0$. As $m$ and $v$ were arbitrary, we conclude
that $\Ric \ge 0$. 

There is a version of the angular Bishop-Gromov inequality for
measured length spaces,
called the ``measure contracting property'' (MCP) \cite{Ohta,Sturm3}.
It satisfies Condition 1. of Wishlist \ref{wishlist}. 

The reason that the MCP notion is not entirely satisfactory can be
seen by asking what it takes
for a smooth measured length space $(M, g, e^{-\Psi}
\: \dvol_M)$ to satisfy the $N$-dimensional angular Bishop-Gromov inequality.
(Here $\dim(M) = n$.)
There is a Riccati-type inequality
\begin{equation} \label{eqn27}
\frac{\partial}{\partial r} \left( \Tr \Pi \: - \: \frac{\partial 
\Psi}{\partial r} \right) \: \le \: - \: \Ric_N(\partial_r, \partial_r) 
\: - \: \frac{1}{N-1} \: \left( \Tr \Pi \: - \: \frac{\partial 
\Psi}{\partial r} \right)^2,
\end{equation}
which looks good.
Again there is an expansion for
the measure of the infinitesimally small angular sector considered above,
of the form 
$\widehat{V}(r) \: = \: r^{n} \left( a_0 \: + a_1 \: r \: + \: a_2 r^2 \: + \:
\ldots \right)$, where the coefficents $a_i$ can be expressed in terms
of curvature derivatives and the derivatives of $\Psi$.
However, if $N > n$ then 
saying that $r^{-N} \widehat{V}(r)$ is nonincreasing in $r$ does not
imply anything about the coefficients. Thus having the $N$-dimensional
angular Bishop-Gromov inequality does not imply that
$\Ric_N \ge 0$.
In particular, it does not seem that one can
prove Corollary \ref{maincor} using MCP. 

Having nonnegative $N$-Ricci
curvature does imply MCP \cite{Sturm3}.
\end{remark}

\section{$N$-Ricci curvature bounded below by $K$} \label{sec5}

In Section \ref{sec3} we gave the definition of nonnegative $N$-Ricci
curvature.  In this section we discuss how to extend this to
a notion of a measured length space having $N$-Ricci curvature
bounded below by some real number $K$.

We start with the case $N=\infty$. As mentioned in Section \ref{sec2},
formal computations indicate that in the case of a smooth measured length
space $(M,g,e^{-\Psi} \: \dvol_M)$, 
having $\Ric_\infty \: \ge \: K g$ should imply that
$H_\infty$ has Hessian bounded below by $K g_{H^{-1}}$ on $P(M)$.
In particular, if $\{\mu_t\}_{t \in [0,1]}$ is a geodesic in 
$P(M)$ then we would
expect that $H_\infty(\mu_t) \: - \: \frac{K}{2} \:
W_2(\mu_0, \mu_1)^2 \: t^2 $ is
convex in $t$. This motivates an adaption of Definition \ref{def3}.

In order to handle all $U \in \DC_\infty$,  we first make the
following definition. Given a continuous convex function $U \:
: \: [0, \infty) \rightarrow \R$, we define its ``pressure'' by
\begin{equation} \label{pdef}
p(r) = r U_+^\prime(r) - U(r),
\end{equation}
where $U_+^\prime(r)$ is the right-derivative.  Then given 
$K \in \R$, we define 
$\lambda \: : \: \DC_\infty \rightarrow \R \cup \{-\infty\}$ by
\begeq\label{lambdainfty} 
\lambda(U) = \inf_{r > 0} K \: \frac{p(r)}r = 
\begin{cases}
K \lim_{r \rightarrow 0^+} \frac{p(r)}r & \text{if $K > 0$}, \\
0 & \text{if $K = 0$}, \\
K \lim_{r \rightarrow \infty} \frac{p(r)}r & \text{if $K < 0$}.
\end{cases}
\endeq
Note that if $U = U_\infty$ (recall that $U_\infty(r) \: = \:
r \: \log r$) then $p(r) = r$ and so
$\lambda(U_\infty) \: = \: K$.

\begin{definition} \cite[Definition 5.13]{LV} \label{def5}
Given $K \in \R$, we say that $(X,d,\nu)$ has $\infty$-Ricci curvature
bounded below by $K$ if for all
$\mu_0, \mu_1 \in P(X)$ with
$\supp(\mu_0) \subset \supp(\nu)$ and $\supp(\mu_1) \subset \supp(\nu)$,
there is {\em some} Wasserstein geodesic
$\{\mu_t\}_{t \in [0,1]}$ from $\mu_0$ to $\mu_1$ so that
for all $U \in \DC_\infty$ and all $t \in [0,1]$,
\begin{equation} \label{ineq2}
U_\nu(\mu_t) \: \le \: t \: U_\nu(\mu_1) \: + \: (1-t) \: 
U_\nu(\mu_0)
 \: - \:
\frac12 \: \lambda(U) \: t(1-t) W_2(\mu_0, \mu_1)^2.
\end{equation}
\end{definition}

\begin{remark} \label{rmk9}
A similar definition, but in terms
of $U = U_\infty$ instead of $U \in \DC_\infty$, was used in
\cite[Definition 4.5]{Sturm2}.
\end{remark}

Clearly if $K = 0$ then we recover the notion of
nonnegative $\infty$-Ricci curvature in the sense of Definition \ref{def3}.
The $N=\infty$ results of Sections
\ref{sec3} and \ref{sec4} 
can be extended to the present case where $K$ may be nonzero.

A good notion of $(X,d,\nu)$ having $N$-Ricci curvature bounded below
by $K \in \R$, where $N$ can be finite, is less clear and is essentially
due to Sturm \cite{Sturm3}. The following definition is a variation
of Sturm's definition and appears in \cite{LV2}.

Given $K\in\R$ and $N \in (1, \infty]$, define
\begin{equation} \label{beta}
\beta_t(x_0,x_1) \: = \: 
\begin{cases} 
e^{\frac16 K\: (1-t^2) \: d(x_0,x_1)^2} \qquad\qquad\qquad & \text{if 
$N = \infty$},\\
\infty \qquad\qquad\qquad & \text{if $N < \infty$, $K>0$ and $\alpha>\pi$},\\
\left(\frac{\sin(t\alpha)}{t\sin\alpha}\right)^{N-1}\qquad
& \text{if $N < \infty$, $K>0$ and $\alpha\in[0,\pi]$}, \\
1 \qquad \qquad\qquad & \text{if $N < \infty$ and $K=0$},\\
\left(\frac{\sinh(t\alpha)}{t\sinh\alpha}\right)^{N-1}\qquad
& \text{if $N < \infty$ and $K<0$},
\end{cases}
\end{equation}
where
\begeq
\alpha \: = \: \sqrt{\frac{|K|}{N-1}}\, d(x_0,x_1).
\endeq
When $N = 1$, define
\begin{equation}
\beta_t(x_0,x_1) \: = \: 
\begin{cases} 
\infty \qquad\qquad\qquad & \text{if 
$K>0$},\\
1\qquad
& \text{if $K \le 0$},
\end{cases}
\end{equation}
Although we may not write it explicitly, $\alpha$ and $\beta$ depend
on $K$ and $N$.

We can disintegrate a transference plan $\pi$ with respect to its
first marginal $\mu_0$ or its second marginal $\mu_1$. We write
this in a slightly informal way:
\begin{equation}
d\pi(x_0, x_1) = d\pi(x_1|x_0) d\mu_0(x_0) = d\pi(x_0|x_1) d\mu_1(x_1).
\end{equation}

\begin{definition} \cite{LV2} \label{def6}
We say that $(X, d, \nu)$ has $N$-Ricci curvature bounded below by $K$
if the following condition is satisfied. Given $\mu_0, \mu_1 \in P(X)$
with support in $\supp(\nu)$, write
their Lebesgue decompositions with respect to $\nu$ as
$\mu_0=\rho_0\,\nu + \mu_{0,s}$ and $\mu_1=\rho_1\,\nu + \mu_{1,s}$,
respectively.
Then there is {\em some} 
optimal dynamical transference plan $\Pi$ from $\mu_0$ to $\mu_1$,
with corresponding Wasserstein geodesic
$\mu_t \: = \: (e_t)_*\Pi$, so that for all $U\in\DC_N$ and all $t\in [0,1]$,
we have
\begin{align} \label{ineqDCN}
U_\nu(\mu_t) \: \leq \: & (1-t) \int_{X \times X} \beta_{1-t}(x_0,x_1) \,
U \left( \frac{\rho_0(x_0)}{\beta_{1-t}(x_0,x_1)} \right)
\, d\pi(x_1|x_0)\,d\nu(x_0) \: +  \\
& t \int_{X \times X} \beta_t(x_0,x_1)\, U \left( \frac{\rho_1(x_1)}
{\beta_{t}(x_0,x_1)} \right) \, d\pi(x_0|x_1)\,d\nu(x_1) \: + \notag \\
& U'(\infty) \bigl[ (1-t) \mu_{0,s}[X] \: + \: t \mu_{1,s}[X] \bigr]. \notag
\end{align}
\end{definition}

Here if $\beta_t(x_0, x_1) = \infty$ then we interpret
$\beta_t(x_0,x_1)\, U \left( \frac{\rho_1(x_1)}
{\beta_{t}(x_0,x_1)} \right)$ as
$U^\prime(0) \: \rho_1(x_1)$, and similarly for
$\beta_{1-t}(x_0,x_1) \,
U \left( \frac{\rho_0(x_0)}{\beta_{1-t}(x_0,x_1)} \right)$.

\begin{remark} \label{rmk10}
If $\mu_0$ and $\mu_1$ are absolutely continuous
with respect to $\nu$ then
the inequality can be rewritten in the more symmetric
form
\begin{align} \label{symform}
U_\nu(\mu_t) \: \leq \: & (1-t) \int_{X \times X} 
\frac{\beta_{1-t}(x_0,x_1)}{\rho_0(x_0)} \,
U \left(\frac{\rho_0(x_0)}{\beta_{1-t}(x_0,x_1)} \right)\,d\pi(x_0, x_1)
\: + \\
& t \int_{X \times X} \frac{\beta_t(x_0,x_1)}{\rho_1(x_1)}\, 
U \left( \frac{\rho_1(x_1)}
{\beta_{t}(x_0,x_1)} \right) \, d\pi(x_0, x_1). \notag
\end{align}
\end{remark}

\begin{remark} \label{rmk11}
Given $K \ge K^\prime$ and $N \le N^\prime$, if
$(X, d, \nu)$ has $N$-Ricci curvature bounded below by $K$ then
it also has $N^\prime$-Ricci curvature bounded below by $K^\prime$.
\end{remark}

\begin{remark} \label{rmk12}
The case $N = \infty$ of Definition \ref{def6} is not quite
the same as what we gave in Definition \ref{def5}!  However, it
is true that having $\infty$-Ricci curvature bounded below by 
$K$ in the sense of Definition \ref{def6} implies that one has
$\infty$-Ricci curvature bounded below by 
$K$ in the sense of Definition \ref{def5} \cite{LV2}. Hence
any $N =\infty$ consequences of Definition \ref{def5} are also
consequences of Definition \ref{def6}.
We include the $N=\infty$ case in Definition \ref{def6}
in order to present a unified treatment, but this example shows
that there may be some flexibility in the precise definitions.
\end{remark}

The results of Sections \ref{sec3} and \ref{sec4} now have extensions
to the case $K \neq 0$. However, the proofs of some of the extensions, 
such as that of Theorem \ref{thm1}, may become much more involved
\cite[Theorem 3.1]{Sturm3},\cite{TOT2}. 

Using the extension of Proposition \ref{prop7}, one obtains a
generalized Bonnet-Myers theorem.

\begin{proposition} \cite[Corollary 2.6]{Sturm3} \label{prop11}
If $(X,d,\nu)$ has $N$-Ricci curvature bounded below by $K > 0$ then
$\supp(\nu)$ has diameter bounded above by
$\sqrt{\frac{N-1}{K}} \: \pi$.
\end{proposition}

\section{Analytic consequences} \label{sec6}

Lower Ricci curvature bounds on Riemannian manifolds have
various analytic implications, such as eigenvalue 
inequalities, Sobolev inequalities and local Poincar\'e
inequalities. It turns out that these inequalities pass
to our generalized setting.

\subsection{Log Sobolev and Poincar\'e inequalities}

Let us first discuss the so-called log Sobolev inequality.
If a smooth measured length space $(M, g, 
e^{- \Psi} \: \dvol_M)$ has $\Ric_\infty \: \ge \:
Kg$, with $K > 0$, then for all 
$f \in C^\infty(M)$ with $\int_M f^2 \: e^{- \Psi} \: \dvol_M
\: = \: 1$, it was shown in \cite{Bakry-Emery (1985)} that
\begin{equation} \label{eqn28}
\int_M f^2 \: \log(f^2)  \: e^{- \Psi} \: \dvol_M \: \le \:
\frac{2}{K} \int_M |\nabla f|^2 \: e^{- \Psi} \: \dvol_M.
\end{equation}
The standard log Sobolev inequality on $\R^n$ comes from taking
$d\nu \: = \: (4\pi)^{- \frac{n}{2}} \: e^{- |x|^2} \: d^nx$,
giving
\begin{equation} \label{eqn29}
\int_{\R^n} f^2 \: \log(f^2)  \: e^{- |x|^2} \: d^n x \: \le \:
\int_{\R^n} |\nabla f|^2 \: e^{- |x|^2} \: d^n x
\end{equation}
whenever $(4\pi)^{- \frac{n}{2}} \: \int_{\R^n} f^2 \:
e^{- |x|^2} \: d^nx \: = \: 1$.

The log Sobolev inequality for $(M, g, 
e^{- \Psi} \: \dvol_M)$ was given both heuristic and
rigorous optimal transport proofs by Otto and Villani
\cite{Otto-Villani (2000)}. We describe the heuristic proof here.
From Section \ref{sec2}, having $\Ric_\infty \: \ge \: Kg$
formally implies that $\Hess(H_\infty) \: \ge \: K
g_{H^{-1}}$ on $P(M)$.
Take $\mu_0 \: = \: e^{- \Psi} \: \dvol_M$ and 
$\mu_1 \: = \: f^2 \: e^{- \Psi} \: \dvol_M$. Let
$\{ \mu_t \}_{t \in [0,1]}$ be a Wasserstein geodesic from
$\mu_0$ to $\mu_1$ along which 
\begin{equation} \label{eqn30}
F(t) \: = \: H_\infty(\mu_t) \: - \:
\frac{K}{2} \: W_2(\mu_0, \mu_1)^2 \: t^2
\end{equation}
is convex in $t$. As $F(0) \: = \: 0$, we have
$F(1) \: \le \: F^\prime(1)$, or
\begin{align}
H_\infty(\mu_1) \: - \: \frac{K}{2} \: W_2(\mu_0, \mu_1)^2 \: & \le \:
\left\langle \frac{d\mu_t}{dt} \Big|_{t=1}, (\grad H_\infty)(\mu_1)
\right\rangle_{g_{H^{-1}}} \: - \: K \: W_2(\mu_0, \mu_1)^2 \\
& \le \:
\Bigg| \frac{d\mu_t}{dt} \Big|_{t=1} \Bigg| \cdot
\Bigg| (\grad H_\infty)(\mu_1) \Bigg| \: 
- \: K \: W_2(\mu_0, \mu_1)^2. \notag
\end{align}
Here $\grad H_\infty$ is the formal gradient of $H_\infty$ on
$P(M)$ and the last norms denote lengths with respect to $g_{H^{-1}}$. 
As $\{ \mu_t \}_{t \in [0,1]}$ is a minimizing geodesic from
$\mu_0$ to $\mu_1$, we should have
\begin{equation} \label{eqn31}
\Bigg| \frac{d\mu_t}{dt} \Big|_{t=1} \Bigg| \: = \:
W_2(\mu_0, \mu_1).
\end{equation} 
A formal computation gives
\begin{equation} \label{eqn32}
\Bigg| (\grad H_\infty)(\mu_1) \Bigg|^2 \: = \: 4 \: 
\int_M |\nabla f|^2 \: e^{- \Psi} \: \dvol_M.
\end{equation}
Then
\begin{align}
\int_M f^2 \: \log(f^2)  \: e^{- \Psi} \: \dvol_M \: & \le \:
2 \: W_2(\mu_0, \mu_1) \: 
\sqrt{\int_M |\nabla f|^2 \: e^{- \Psi} \: \dvol_M} \: - \:
\frac{K}{2} \: W_2(\mu_0, \mu_1)^2 \\
& \le \:
\sup_{w \in \R} \left( 
2 \: w \: 
\sqrt{\int_M |\nabla f|^2 \: e^{- \Psi} \: \dvol_M} \: - \:
\frac{K}{2} \: w^2 \right) \notag \\
&  = \: 
\frac{2}{K} \int_M |\nabla f|^2 \: e^{- \Psi} \: \dvol_M \notag
\end{align}
which is the log Sobolev inequality.

The rigorous optimal transport proof in \cite{Otto-Villani (2000)} 
extends to
measured length spaces. To give the statement,
we first must say what we mean by
$|\nabla f|$.
We define the {\em local gradient norm} of a Lipschitz function 
$f \in \Lip(X)$ by the formula
\begin{equation} \label{eqn33}
|\nabla f| (x) \: = \:
\limsup_{y\to x} \frac{|f(y)-f(x)|}{d(x,y)}.
\end{equation}
We don't claim to know the meaning of the gradient $\nabla f$ on $X$ in
this generality, but we can talk about its norm anyway!
Then we have the following log Sobolev inequality for
measured length spaces.

\begin{theorem} \cite[Corollary 6.12]{LV} \label{thm3}
Suppose that a compact measured length space
$(X,d,\nu)$ has $\infty$-Ricci curvature bounded below by
$K > 0$, in the sense of Definition \ref{def5}. Suppose that
$f \in \Lip(X)$ satisfies $\int_X f^2 \: d\nu \: = \: 1$.
Then
\begin{equation} \label{introeqn}
\int_X f^2 \: \log(f^2) \: d\nu \: \le \: \frac{2}{K} \:
\int_X |\nabla f|^2 \: d\nu.
\end{equation}
\end{theorem}

In the case of Riemannian manifolds, one recovers from 
(\ref{introeqn}) the log Sobolev inequality (\ref{eqn28})
of Bakry and \'Emery.

\begin{remark} \label{rmk13} 
The proof of Theorem \ref{thm3}, along with the other
inequalities in this section, uses the $K>0$ analog of Proposition \ref{prop8}.
In turn,
the proof of Proposition \ref{prop8} uses the fact that (\ref{ineq1})
holds for all $U \in \DC_N$, as opposed to just $U_N$.
\end{remark}

As is well-known, one can obtain a Poincar\'e inequality
from (\ref{introeqn}). Take $h \in \Lip(X)$ with $\int_X h\,d\nu = 0$
and put $f^2 \: = \: 1 + \epsilon h$. Taking
$\epsilon$ small and expanding the two sides of
(\ref{introeqn}) in $\epsilon$ gives the following result.

\begin{corollary} \cite[Theorem 6.18]{LV} \label{corr}
Suppose that a compact measured length space
$(X,d,\nu)$ has $\infty$-Ricci curvature bounded below by
$K > 0$. Then for all
$h \in \Lip(X)$ with $\int_X h\,d\nu = 0$, we have
\begin{equation} \label{eqn34}
\int_X h^2\,d\nu \leq 
\frac1{K} \int_X |\nabla h|^2\,d\nu. 
\end{equation}
\end{corollary}

In case of a smooth measured length space
$(M, g, e^{-\Psi} \: \dvol_M)$, the inequality (\ref{eqn34}) 
coincides with
the Bakry-\'Emery extension of the Lichnerowicz inequality,
namely $\lambda_1 (\widetilde{\triangle}) \: \ge \: K$.
For a general measured length space as in the hypotheses of
Corollary \ref{corr}, we do not know if there is a well-defined
Laplacian.  The Poincar\'e inequality of Corollary \ref{corr}
can be seen as a generalized eigenvalue inequality that
avoids this issue.  To say a bit more about when one does
have a Laplacian,
if $Q(h) \: = \: \int_X |\nabla h|^2 \: d\nu$ defines a quadratic form on
$\Lip(X)$, which in addition is closable in $L^2(X, \nu)$, then
there is a self-adjoint Laplacian $\triangle_\nu$
associated to $Q$. In this case, Corollary \ref{corr} implies that 
$\triangle_\nu \: \ge \: K$ on the orthogonal complement of the constant
functions.

In the case of a Ricci limit space, Cheeger and Colding used 
additional structure in order to show
the Laplacian does exist \cite{Cheeger-Colding (2000)}.

\subsection{Sobolev inequality}

The log Sobolev inequality can be viewed as an 
infinite-dimensional version of an ordinary Sobolev inequality.
As such, it is interesting because it is a dimension-independent
result.  However, if one has $N$-Ricci curvature 
bounded below by $K > 0$ with $N$ finite then one gets an
ordinary Sobolev inequality, which is a sharper result.

\begin{proposition} \cite{LV2} \label{prop12}
Given $N \in (1, \infty)$ and $K > 0$, suppose that $(X,d,\nu)$ has $N$-Ricci
curvature bounded below by $K$.
Then for any nonnegative Lipschitz function $\rho_0 \in \Lip(X)$ with
$\int_X \rho_0 \: d\nu \: = \: 1$, one has
\begin{equation} \label{Sobolev}
N - N \int_X \rho_0^{1-\frac{1}{N}} \: d\nu \: \le \:
\frac{1}{2K} \left( \frac{N-1}{N} \right)^2 \int_X
\frac{\rho_0^{-1-\frac{2}{N}}}{\frac13 + \frac23 \rho_0^{-\frac{1}{N}}} \:
|\nabla \rho_0|^2 \: d\nu.
\end{equation}
\end{proposition}

To put Proposition \ref{prop12} into a more conventional form,
we give a slightly weaker inequality.

\begin{proposition} \cite{LV2} \label{prop13}
Given $N \in (2, \infty)$ and $K > 0$, suppose that $(X,d,\nu)$ has $N$-Ricci
curvature bounded below by $K$.
Then for any nonnegative Lipschitz function $f \in \Lip(X)$ with
$\int_X f^{\frac{2N}{N-2}} \: d\nu \: = \: 1$, one has
\begin{equation} \label{Sobolev2}
1 -  \left( \int_X f \: d\nu \right)^{\frac{2}{N+2}} \: \le \:
\frac{6}{KN} \left( \frac{N-1}{N-2} \right)^2 \int_X
|\nabla f|^2 \: d\nu.
\end{equation}
\end{proposition}

Putting (\ref{Sobolev2}) into a homogeneous form, 
the content of Proposition \ref{prop13} is that there
is a bound of the form
$\parallel f \parallel_{\frac{2N}{N-2}} \: \le
F \left( \parallel f \parallel_1,
\parallel \nabla f \parallel_2 \right)$ for some
appropriate function $F$. This is
an example of Sobolev embedding. The inequality
(\ref{Sobolev2}) is not sharp, due to the
many approximations made in its derivation.

One can use Proposition \ref{prop12} to prove 
a sharp Poincar\'e inequality.

\begin{proposition} \cite{LV2} \label{prop14}
Given $N \in  (1, \infty)$ and $K > 0$, 
suppose that $(X,d, \nu)$ has $N$-Ricci curvature
bounded below by $K$. Suppose that $h \in \Lip(X)$ has
$\int_X h \: d\nu \: = \: 0$. Then
\begin{equation} \label{eqn35}
\int_X h^2\,d\nu \leq 
\frac{N-1}{KN} \int_X |\nabla h|^2\,d\nu. 
\end{equation}
\end{proposition}

In the case of an $N$-dimensional Riemannian manifold with 
$\Ric \: \ge \: K \: g$,
one recovers the Lichnerowicz inequality for the lowest positive
eigenvalue of the Laplacian
\cite{Lichnerowicz}. It is sharp on round spheres.

\subsection{Local Poincar\'e inequality}

When doing analysis on metric-measure spaces, a useful analytic property
is a ``local'' Poincar\'e inequality.
A metric-measure space $(X, d, \nu)$ admits a local Poincar\'e inequality 
if, roughly speaking, for each function
$f$ and each ball $B$ in $X$, the mean deviation (on $B$) of $f$
from its average value on $B$ is quantitatively controlled by the 
gradient of $f$ on a larger ball.

To make this precise, if $B = B_r (x)$ is a ball in $X$ then we write
$\lambda B$ for $B_{\lambda r}(x)$.
The measure $\nu$ is said to be {\em doubling} if
there is some $D > 0$ so that for all balls $B$,
$\nu(2B) \: \le \: D \: \nu(B)$. 
An {\em upper gradient} for a function
$u \in C(X)$ is a Borel function $g \: : \: X \rightarrow
[0, \infty]$ such that for each curve
$\gamma \: : \: [0,1] \rightarrow X$ with finite length 
$L(\gamma)$ 
and constant speed,
\begin{equation} \label{eqn36}
\bigl| u(\gamma(1)) \: - \: u(\gamma(0)) \bigr| \: \le \:
L(\gamma) \: \int_0^1 \: g(\gamma(t)) \: dt.
\end{equation}
If $u$ is Lipschitz then $|\nabla u|$ is an example of an upper gradient.

There are many forms of local Poincar\'e inequalities. 
The strongest one, in a certain sense, is as follows :
\begin{definition} \label{def7}
A metric-measure space $(X, d, \nu)$ admits a local Poincar\'e
inequality if
there are constants $\lambda\geq 1$ and 
$P<\infty$ such that for all $u\in C(X)$ and $B=B_r(x)$ with
$\nu(B) > 0$,
each upper gradient $g$ of $u$ satisfies
\begeq\label{PI2}
\barint_B |u-\av{u}_B|\,d\nu \leq P r\, \barint_{\lambda B}
g \,d\nu.
\endeq
\end{definition}

Here the barred integral is the average (with respect to $\nu$),
e.g. $\barint_{\lambda B}
g \,d\nu \: = \: \frac{\int_{\lambda B}
g \,d\nu}{\nu(\lambda B)}$,  and
$\av{u}_B$ is the average of $u$ over the ball $B$. In the case of
a length space, the local
Poincar\'e inequality as formulated in Definition \ref{def7}
actually implies stronger inequalities, for which we refer to
\cite[Chapters 4 and 9]{Heinonen}.  It is known that the property
of admitting a local Poincar\'e inequality is preserved under
measured Gromov-Hausdorff limits \cite{Keith,Koskela}.
(This was also shown by Cheeger in unpublished work.)
Cheeger showed that if a metric-measure space
has a doubling measure and admits a local Poincar\'e inequality then
it has remarkable extra local structure 
\cite{Cheeger}.

Cheeger and Colding showed that local Poincar\'e inequalities exist for
Ricci limit spaces \cite{Cheeger-Colding (2000)}.  The method of proof
was to show that Riemannian manifolds with lower Ricci curvature
bounds satisfy a certain ``segment
inequality''
\cite[Theorem 2.11]{Cheeger-Colding}
and then to show that the property of satisfying
the segment inequality is preserved under measured Gromov-Hausdorff 
limits \cite[Theorem 2.6]{Cheeger-Colding (2000)}. 
The segment inequality then implies
the local Poincar\'e inequality.

It turns out that the argument
using the segment inequality can be abstracted and
applied to certain measured length spaces. For simplicity, we
restrict to the case of nonnegative $N$-Ricci curvature.
We say that $(X, d, \nu)$ has
almost-everywhere unique geodesics if
for $\nu\otimes\nu$-almost all $(x_0,x_1)\in X\times X$,
there is a unique minimizing geodesic
$\gamma \in \Gamma$ with $\gamma(0)=x_0$ and
$\gamma(1)=x_1$.

\begin{theorem} \cite{LV2,vR,Sturm3} \label{thm4}
If a compact measured length space 
$(X,d,\nu)$ has nonnegative $N$-Ricci curvature
and almost-everywhere unique geodesics then it satisfies
the local Poincar\'e inequality of Definition~\ref{def7} with
$\lambda\: = \: 2$ and $P\: = \: 2^{2N+1}$.
\end{theorem}

As is well-known, a Riemannian manifold has almost-everywhere
unique geodesics. A sufficient condition for $(X,d, \nu)$ to
have almost-everywhere unique geodesics is that 
almost every $x \in X$ is nonbranching in a certain sense
\cite{vR,Sturm3}.

The result of Theorem \ref{thm4} holds in greater generality.  What one needs
is a way of joining up points by geodesics, called a ``democratic
coupling'' in \cite{LV2}, and a doubling condition on the measure.

We do not know whether the condition of 
nonnegative $N$-Ricci curvature is sufficient in itself to imply a
local Poincar\'e inequality. Having nonnegative $N$-Ricci
curvature does not imply almost-everywhere unique geodesics.
For a noncompact example, the finite-dimensional Banach space $\R^n$ with
the $l_1$ norm and the Lebesgue measure has nonnegative
$n$-Ricci curvature, but certainly does not have almost-everywhere
unique geodesics.

\section{Final remarks} \label{sec7}
In this survey we have concentrated on compact spaces.
There is also a notion of Ricci curvature bounded below
for noncompact measured length spaces $(X,d,\nu)$
\cite[Appendix E]{LV}. Here we want $X$ to be a complete pointed locally
compact length space and $\nu$ to be a nonnegative nonzero 
Radon measure on $X$. We do not require $\nu$ to be a
probability measure.  There is a Wasserstein space $P_2(X)$
of probability measures on $X$ with finite second moment,
i.e. 
\begin{equation} \label{eqn37}
P_2(X) \: = \: \left\{ \mu \in P(X) \: : \: \int_X d(\bp,x)^2 \: d\mu(x) 
\: < \: \infty \right\},
\end{equation}
where $\bp$ is the basepoint in $X$. Many of the results described
in this survey extend from compact spaces to such noncompact spaces,
although interesting technical points arise.

In particular, if $(X, d, \nu)$ is a compact or
noncompact space with nonnegative $N$-Ricci
curvature and $\supp(\nu) = X$, and if $x$ is a point in $X$, then
a tangent cone at $x$ has nonnegative $N$-Ricci curvature
\cite[Corollary E.44]{LV}.

There are many directions for future research. Any
specific problems that we write here may become obsolete, but let us just
mention two general directions.  One direction is to
see whether known results about Riemannian manifolds with
lower Ricci curvature bounds extend to measured length spaces
with lower Ricci curvature bounds.  As a caution, this is not
always the case. For example, the Cheeger-Gromoll splitting
theorem says that if there is a line in a complete Riemannian manifold 
$M$ with nonnegative Ricci curvature then there is an 
isometric splitting $M = \R \times Y$. This is not true for
measured length spaces with nonnegative $N$-Ricci curvature.
Counterexamples are given by nonEuclidean $n$-dimensional normed
linear spaces,
equipped with Lebesgue measure, which all have nonnegative
$n$-Ricci curvature \cite{TOT2}. However, it is possible that
there is some vestige of the splitting theorem left.

The splitting theorem does hold for a pointed Gromov-Hausdorff limit
of a sequence $\{(M_i, g_i)\}_{i=1}^\infty$ of complete Riemannian
manifolds with Ricci curvature bounded below by $- \: \frac{1}{i}$
\cite{Cheeger-Colding}, so not every finite-dimensional
$(X,d,\nu)$ with nonnegative
$N$-Ricci curvature arises as a limit in this way.
(The analogous statement is not known for finite-dimensional
Alexandrov spaces,
but there are candidate Alexandrov spaces that may not
be Gromov-Hausdorff limits of Riemannian manifolds with
sectional curvature uniformly bounded below \cite{Kapovitch}.)
One's attitude towards this fact may depend on whether one intuitively
feels that finite-dimensional normed linear spaces should 
or should not have nonnegative Ricci curvature.

Another direction of research is to find classes of measured length
spaces $(X,d,\nu)$ which do or do not have lower Ricci curvature bounds.
This usually amounts to understanding optimal transport on
such spaces.

\bibliographystyle{acm}

\end{document}